\documentstyle[amssymb,amscd,12pt]{amsart}

\newtheorem{defn}{Definition}
\newtheorem{thm}[defn]{Theorem}

\newtheorem{cor}[defn]{Corollary}
\newtheorem{lem}[defn]{Lemma}

\theoremstyle{remark}
\newtheorem{rem}{Remark}
\theoremstyle{remark}
\newtheorem{exam}{Example}

\numberwithin{equation}{section}
\numberwithin{defn}{section}

\begin{document}


\newcommand\spk{{\operatorname{Spec}}(k)}
\renewcommand\sp{\operatorname{Spec}}
\renewcommand\sf{\operatorname{Spf}}
\newcommand\proj{\operatorname{Proj}}
\newcommand\aut{\operatorname{Aut}}
\newcommand\grv{{\operatorname{Gr}}(V)}
\newcommand\gr{\operatorname{Gr}}
\newcommand\glv{{\operatorname{Gl}}(V)}
\newcommand\glve{{\widetilde{\operatorname{Gl}}(V)}}
\newcommand\gl{\operatorname{Gl}}
\renewcommand\hom{\operatorname{Hom}}
\renewcommand\det{\operatorname{Det}}
\newcommand\detd{\operatorname{Det}^\ast}
\newcommand\im{\operatorname{Im}}
\newcommand\Pf{\operatorname{Pf}}
\newcommand\res{\operatorname{Res}}
\newcommand\pic{\operatorname{Pic}}
\newcommand\limi{\varinjlim}
\newcommand\limil[1]{\underset{#1}\varinjlim\,}
\newcommand\limp{\varprojlim}
\newcommand\limpl[1]{\underset{#1}\varprojlim\,}

\renewcommand\o{{\cal O}}
\renewcommand\L{{\cal L}}
\renewcommand\c{{\cal C}^\bullet}
\renewcommand\P{{\cal P}}
\newcommand\Z{{\Bbb Z}}
\newcommand\A{{\Bbb A}}
\newcommand\M{{\cal M}}

\renewcommand\tilde{\widetilde}

\newcommand\cur{{\cal M}_{\infty}}
\newcommand\pry{\tilde{\P}_{\infty}}

\newcommand\iso{@>{\sim}>>}
\renewcommand\lim{\varprojlim \Sb A \\ A\sim V_+ \endSb}
\newcommand\fu{\underline}
\newcommand\kz{\fu{k((z))}^*}
\newcommand\beq{
      \setcounter{equation}{\value{defn}}\addtocounter{defn}1
      \begin{equation}}

\renewcommand{\thesubsection}{\thesection.\Alph{subsection}}
\renewcommand{\thefootnote}{\alph{footnote}}

\hfill{\fbox{\tiny To appear in {\bf  International  Journal of
Mathematics}}}
\vskip1truecm

\title [Prym Varieties and the Infinite Grassmannian]
{Prym Varieties and the \\ Infinite Grassmannian}

\author[F. J. Plaza Mart\'{\i}n] {Francisco J. Plaza
Mart\'{\i}n\\ \medskip\tiny Departamento de Matem\'atica Pura y
Aplicada\\ Universidad de Salamanca}

\address{ Departamento de Matem\'atica Pura y Aplicada \\ Universidad de
Salamanca \\ Plaza de la Merced 1-4\\37008 Salamanca. Spain.}

\thanks{1980 Mathematics Subject Classification (1985 Revision).
Primary 14H40. Secondary 35Q20. \\
This work is partially supported by the
CICYT research contract n. PB91-0188.}

\email{fplaza@@gugu.usal.es}


\begin{abstract}
In this paper we study Prym varieties and their moduli space using the
well known techniques of the infinite Grassmannian. There are three main
results of this paper: a new definition of the BKP hierarchy over an
arbitrary base field (that generalizes the classical one over
${\mathbb C}$);  a characterization of Prym varieties in terms of
dynamical systems, and  explicit equations for the moduli space of
(certain) Prym varieties. For all of these problems the language of the
infinite Grassmannian, in its algebro-geometric version, allows us to
deal with these problems from the same point of view.
\end{abstract}

\maketitle


\setcounter{tocdepth}1
\tableofcontents

\section{Introduction}

The aim of this paper is to generalize some results concerning the BKP
hierarchy and geometric characterizations of Jacobians and Pryms
proved in \cite{LiMu,Mul1,Sh,Sh2} and to study the moduli space of
(certain) Prym varieties following similar ideas to those of \cite{MP}.
I should remark that the techniques employed here are those of
algebraic geometry, and most statements are therefore valid over an
arbitrary base field. The organization of the paper is as follows:

In \S2 some basic definitions, results and tools needed for the next
sections are introduced. Some of them are known (e.g. infinite
Grassmannian, Krichever functor, $\Gamma$ group, etc) and have their
origin in the study of the moduli space of Riemann surfaces, Jacobian
varieties and conformal field theory (\cite{BS,DKJM,KSU,Mul1,N,PS,SW} .
To develop analogues for the theory of Prym varieties, we
will define certain subschemes of the infinite Grassmannian and
will avoid the introduction of the formalism of the n-component KP
hierarchy. However, the analogue of the determinant bundle is a
non-trivial problem. An important result ({\ref{thm:pfaffian}}) of this
section is the existence of a square root of the determinant bundle
(over a certain subscheme), which will be called Pfaffian, as well as
an explicit construction of global sections of this Pfaffian bundle,
improving previous results of \cite{B,PS}.

Section \S3 contains a new definition of the BKP hierarchy as the
defining equations of a suitable subscheme of the infinite
Grassmannian ({\ref{defn:BKP}}); namely, the space of maximal totally
isotropic subspaces.  This new definition is quite natural since we
should recall that the KP hierarchy is in fact equivalent to the
Pl\"ucker equations, which are the defining equations of the infinite
Grassmannian (\cite{SS}). The definition is therefore valid over an
arbitrary base field, but is shown to be equivalent to the
classical one when the base field is
$\mathbb C$. We wish  to point out that this new definition of the BKP
hierarchy, together with the previously defined Pfaffian line bundle,
will give a justification for some results relating
solutions of the KP and BKP hierarchies (\cite{DKJM,DKJM2}) as well as
for techniques involving pfaffians when computing $\tau$-functions for
the BKP hierarchy (\cite{H,O}).  Note that this (algebraic) approach
allows us to introduce the BKP hierarchy with no mention of
pseduodifferential operators.

The standard way to relate the study of Jacobians to that of infinite
Grassmannian is the Krichever functor (\cite{Mul2}). Two remarkable
papers dealing with the case of Pryms are \cite{LiMu,Sh2}. The
characterization given in \cite{LiMu} for a point of (a certain
quotient of) the infinite Grassmannian to be associated to a Prym (via
the Krichever functor) is that its orbit (under a suitable group) is
finite dimensional (see Theorem 5.14 of \cite{LiMu} for the precise
statement). Moreover, that action is interpreted as a dynamical system
on that space. The idea of Shiota (\cite{Sh,Sh2}), related to
the above one, is to characterize Jacobian varieties as compact
solutions of the KP hierarchy through the (analytic)  study of
infinitesimal deformations of sheaves.

The characterizations of Jacobians ({\ref{thm:shiota}}) and Pryms
({\ref{thm:shiota-prym}}) given in \S4 profit from both approaches.
Recalling that the action of the formal Jacobian (formal group
$\Gamma$, \S{\ref{subsec:gamma}}) on the Grassmannian is algebraic and
studying the structure of the orbits ({\ref{lem:orbit}}), one can give
an algebraic statement generalizing Theorem 6 of \cite{Sh} and Theorem
5.14 of \cite{LiMu}. But in our statement no quotient space
is needed. Bearing this in mind, there is no problem in extending this
result to the case of Pryms ({\ref{thm:shiota-prym}}). Nevertheless, we
shall use the language of dynamical systems to express both
characterizations.

In the last section (\S5), and following the spirit of \cite{MP},
explicit equations for the moduli space of Pryms are given; first as
Bilinear Identities ({\ref{thm:bil-ident}}), then as partial
differential equations ({\ref{thm:pde-p0}}) when $char(k)=0$. This set
of equations should not be confused with the BKP hierarchy. Here we
make two considerations; firstly, that a formal trivilization is
attached to each datum (e.g. a curve, a bundle, \dots~) since we wish
to work in a uniform frame (e.g. $\gr(k((z)))$, $\Gamma$,
\dots~) and, secondly, that not all Pryms are considered; only those
coming from an integral curve together with an involution with at least
one fixed smooth point are taken into account since, for technical
reasons, the involution should correspond to an automorphism of
$k((z))$ preserving $k[[z]]$.

We address the reader to \cite{AMP} for a detailed discussion of the
infinite Grassmannian and to \cite{Mum} for the basic facts on Prym
varieties needed.

\section{Preliminaries on the Infinite Grassmannian}

\subsection{Basic Facts}
Recall (\cite{AMP, BS}) that given a pair $(V,V_+)$ consisting of a
$k$-vector space\footnote{For simplicity's sake we will
assume that $k$ is an algebraically closed field} and a subspace of it,
there exists a scheme
$\gr(V,V_+)$ over $\spk$  whose rational points are:
$$\left\{\begin{gathered}
\text{ subspaces $L\subseteq \hat V$, such that $L\cap\hat V_+$ }\\
\text{ and ${\hat V}/{L+\hat V_+}$ are of finite dimension }
\end{gathered}\right\}$$
where $\,\hat{}\,$ denotes the completion with respect to the topology
given by the subspaces that are commensurable with $V_+$. The points of
$\gr(V,V_+)$ will be called discrete subspaces. The essential fact for
its existence is that there is a covering by open subfunctors
$F_A$ (where
$A\sim V_+$ are commensurable) representing the functor
$\fu\hom(L_A,\hat A)$, where  $L_A$ is a rational point of $\gr(V,V_+)$
such that $L_A\oplus
\hat A\simeq \hat V$ (\cite{BS}). (See \cite{AMP} for the definition of
the functor of points of $\gr(V,V_+)$). Let us denote this infinite
Grassmannian simply by
$\grv$, and let $\L$ be the universal discrete submodule of
$\hat V_{\grv}$. We will assume that $V$ is complete with respect to
the $V_+$-topology. From this construction it is easily deduced that
$\grv$ is locally integral and separated.

The connected components of $\grv$ are given by the Euler
characteristic (index) of the complex:
\beq
\L\to (\hat V/\hat V_+)_{\grv}
\label{eq:complex}\end{equation}
and that of index $n$ will be denoted by $\gr^n(V,V_+)$.  It is also
shown, that $\grv$ carries a line bundle, $\det_V$, given by the
determinant of the complex {\ref{eq:complex}} (see \cite{KM} for a
general theory of determinants) whose stalk at a rational point $L$ is:
$$\overset{max}\wedge (L\cap\hat V_+)\otimes\overset{max}\wedge ( {\hat
V}/{(L+\hat V_+)})^\ast$$

This line bundle has no global sections but its dual does. Moreover,
for each $A\sim V_+$ one can define a global section, $\Omega_A$ of
$\detd_V$ such that it vanishes outside $F_A$. For every subspace
$\Omega$ of $H^0(\gr^0(V),\detd_V)$ one has a sheaf homomorphism:
$$\Omega\otimes_k\o_{\grv}\to\detd_V$$
If it is surjective ($\Omega$ is ``big enough''), it induces a scheme
homomorphism:
$${\frak p}_V: \gr^0(V) \to \check{\mathbb P}\Omega^*
\,\overset{\text{\tiny{def}}}{=}\,\operatorname{Proj}S^\bullet\Omega$$
which is known as the Pl\"ucker morphism.

Although $k((z))$ is not complete with respect to the $k[[z]]$-topology,
it is easy to see that the functor:
\beq
S\rightsquigarrow\{ L\in\gr(k((z)),k[[z]])(S)\;\vert\; L\subseteq
\o_S((z))\}
\label{eq:grkz}\end{equation}
is locally closed, and it is therefore representable by a closed
subscheme, which will be denoted by $\gr(k((z)))$ again, when no
confussion arises (\cite{P}). It can now be shown that
$H^0(\gr^0(k((z))),\detd_V)$ has a dense subspace, $\Omega$, consisting
of sections of the kind $\Omega_A$ and labelled by Young diagrams
(see \cite{P,SS}), and that the Pl\"ucker morphism is a closed
immersion. The general results given along this section are valid for
 subscheme {\ref{eq:grkz}}.

\subsection{Important subschemes of the Infinite Grassmannian}

Let us now introduce two closed subschemes of
$\grv$ that will be useful in the next sections. Assume now that an
automorphism
$\sigma$ of $V$ (as $k$-vector space) is given, and that $\sigma(\hat
V_+)=\hat V_+$; it then induces an automorphism of the scheme $\grv$.
Since $\grv$ is separated, one has that:
$$\gr_\sigma(V)=\left\{ L\in\grv\,\vert\, \sigma(L)=L\right\}$$
is a closed subscheme of $\grv$.

For the second, recall the isomorphism
$\gr(V,V_+)\iso\gr(V^*,V_+^\circ)$ given by incidence (see \cite{P});
that is, it sends a discrete subspace $L$ to $L^\circ$, the space of
continuous linear forms that vanish on $L$. Let $p: V\iso V^*$ be the
isomorphism of $V$ with its dual vector space $V^*$ induced by an
irreducible hemisymmetric metric on $V$. Assume further that $p(\hat
V_+)=\hat V_+^\circ$. It then induces an isomorphism
$\gr(V^*,V_+^\circ)\iso\gr(V,V_+)$, which composed with the one given by
incidence, gives rise to the following automorphism of $\grv$:
$$\aligned R:\grv &\longrightarrow \grv \\
L &\longmapsto L^{\perp} \endaligned$$
(where $\perp$ denotes the orthogonal with respect to the metric).

Straightforward calculation shows that $R^*\det_V\simeq \det_V$, and that
the index of a point $L\in\grv$ is exactly the opposite of the index of
$R(L)=L^\perp$.

Given $\sigma\in\aut_{k-alg}k((z))$ such that $\sigma^2=Id$, consider
 the following irreducible hemisymmetric metric:
\beq
\begin{aligned} V\times V&\to k \\
(f,g)&\mapsto\res_{z=0} f(z)\cdot (\sigma^*g(z)) dz\end{aligned}
\label{eq:metric}\end{equation}
It is now clear that there exists a closed subscheme $\gr^I_\sigma(V)$ of
$\gr^0(V)$ such that:
{\small \beq\gr^I_\sigma(V)^\bullet(k)=\left\{ L\in\grv\,\vert\,
\begin{gathered}L\text{ is maximal totally isotropic }\\
\text{with respect to the metric {\ref{eq:metric}}}\end{gathered}
\right\}\label{defn:gri}\end{equation}}
From now on, a subspace of $V$ will be called m.t.i. when it is
maximal totally isotropic (compare with \S2.2 of \cite{Sh2}).

\begin{rem}
Whenever we use a hemisymmetric metric or consider the automorphism of
$k((z))$ induced by $z\mapsto -z$, is assumed $char(k) \neq 2$.
\end{rem}

\begin{exam}\label{exam:metric}
This is a fundamental example because it will be the situation when
studying Prym varieties in terms of $\gr^I_\sigma(V)$.
 Let $V=k((z))$, $V_+=k[[z]]$ and $\sigma_0$ that given by $z\mapsto
-z$. The metric is now: $<f(z),g(z)>=\res_{z=0} f(z)g(-z)dz$. It is then
easy to prove that:
{\small $$\gr_0(k((z)),k[[z]])\simeq\gr(k((z^2)),k[[z^2]])
\underset\spk\times\gr(z\cdot k((z^2)),z \cdot k[[z^2]])$$}
(we simply write $0$ instead of $\sigma_0$). It is worth comparing the 2
component BKP hierarchy (as given in
\S3 of \cite{Sh2}) with $\gr_0(k((z)),k[[z]])$.
\end{exam}

\subsection{Pffafian Line Bundle}

\begin{thm}\label{thm:pfaffian}
There exists a line bundle, $\Pf$ (called Pfaffian), over $\gr^I_\sigma(V)$
such that:
$$\Pf^{\otimes 2}\simeq \detd_V\vert_{\gr^I_\sigma(V)}$$
\end{thm}

\begin{pf}
Note that $\detd_V$ is isomorphic to the determinant of the dual of the
complex over $\gr^I_\sigma(V)$: $\hat V_+\to \hat V/\L$ ($\L$ being the
universal m.t.i. submodule). Let us compute its cohomology. Let
$p:V\to V^*$ be the isomorphism induced by the metric
{\ref{eq:metric}}. Then, in the commutative diagram:
$$\CD
 @. \L @>>> \hat V/\hat V_+ @>>> \frac{\hat V}{\hat V_++\L} @>>> 0 \\
 @. @V{p}VV @V{p}VV @VVV \\
0 \to (\frac{\hat V}{\hat V_++\L})^* @>>> (\hat V/\L)^* @>>> \hat V_+^*
@>>> \frac{\hat V_+^*}{(\hat V/\L)^*}@>>> 0
\endCD$$
the two middle vertical arrows are isomorphisms, since $\L$ and $\hat
V_+$ are m.t.i.~; and thus the right one is an isomorphism.

One now has:
$$\detd_V\iso {\left(\bigwedge \hat V/(\hat V_++\L)\right)^*}^{\otimes
2}$$ By the local structure of $\gr^I_\sigma(V)$ it is not difficult to
show that
$\hat V/(\hat V_++\L$ is locally free of finite type, and by \cite{KM}
it makes sense to define:
$$\Pf\,\overset{\text{\tiny{def}}}{=}\,(\bigwedge \hat V/(\hat
V_++\L))^*$$
\end{pf}

An analytic construction of the Pfaffian line bundle can be found in
\cite{B}, but nevertheless we prefer to continue with the
algebraic machinery. Another construction is given in \cite{PS}.

We are now interested in building sections of this line bundle. First,
observe that the covering $\{F_A\}$ ($A\sim V_+$) of $\grv$ induces
another one of $\gr^I_\sigma(V)$ by open subsets of the form $\{\bar
F_A=F_A\cap\gr^I_\sigma(V)\}$ where $F_A\subset\grv$ ($A\sim V_+$) and
$A$ is m.t.i.~. The second ingredient is the following:

\begin{lem}
Let $A$ and $B$ be two m.t.i. subspaces of $V$. One then has a canonical
isomorphism:
$$B/(B\cap A)\iso (A/(B\cap A))^*$$
induced by the metric.
\end{lem}

\begin{pf}
Note that the morphism $p:V\to  V^*$ gives an isomorphism from
$A$ (resp. $B$) to $A^\circ$ (resp. $B^\circ$). Therefore, from the
diagram:
$$\CD 0 @>>> A\cap B @>>> A @>>> A/(A\cap B) @>>> 0 \\
@. @V{p}VV @V{p}VV @VVV \\
0 @>>> B^\circ @>>> V^* @>>> B^* @>>> 0 \endCD$$
one has an injection $A/(A\cap B) \to B^*$. The linear forms belonging
to the image vanish on $A\cap B$, and hence:
$$A/(A\cap B) \to (A\cap B)^\circ\cap B^*\simeq \left(B/(A\cap
B)\right)^*$$
Analogously, one obtains another injection $B/(A\cap B) \to
\left(A/(A\cap B)\right)^*$, and they are transposed of each other, and
are therefore both isomorphisms.
\end{pf}

\begin{thm}
To each $A\sim V_+$ m.t.i. one associates a section
$\bar\Omega_A$ of $\Pf$ that vanishes outside $\bar F_A$, and hence:
$$\Omega_A\vert_{\gr^I_\sigma}=\lambda\cdot\bar\Omega_A^2\qquad
\lambda\in k^*$$

\end{thm}

\begin{pf}
Observe now that by similar arguments to those of \cite{AMP} one can
construct sections of $\Pf$ for each $A\sim V_+$ m.t.i.~. The only
remarkable aspect is that this is possible because
$B/(B\cap A)$ and $(A/(B\cap A))^*$ are canonically isomorphic (for $A,B$
m.t.i. and commensurable with $V_+$).

From this last property, and from the fact
$\Pf^{\otimes 2}\simeq \detd_V$, the claim follows.
\end{pf}

\subsection{Krichever Functor}

\begin{defn}
Define the moduli functor of pointed curves $\cur$ over the category of
$k$-schemes by the sheafication of:
$$S\rightsquigarrow\{\text{ families $(C,D,z)$ over $S$
}\}/\text{equivalence}$$  where these families satisfy:
\begin{enumerate}
\item $\pi:C\to S$ is a proper flat morphism,
whose geometric fibres are integral curves,
\item $s:S\to C$ is a section of $\pi$, such that when
considered as a Cartier Divisor $D$ over $C$ it is smooth and of relative
degree 1, and flat over
$S$. (We understand that $D\subset C$ is smooth over $S$, iff for
every closed point $x\in D$ there exists an open neighborhood $U$ of
$x$ in $C$ such that the morphism $U\to S$ is smooth).
\item $z$ is a formal trivialization of $C$ along $D$; that is, a family
of epimorphisms  of rings:
$$\o_C\longrightarrow
s_*\left({\o_S[t]}/{t^m\,\o_S[t]}\right)\qquad m\in{\mathbb N}$$
 compatible with respect to the canonical projections:
$${\o_S[t]}/{t^m\,\o_S[t]}\to {\o_S[t]}/{t^{m'}\,\o_S[t]}\qquad m\ge
m'$$
and such that that corresponding to $m=1$ equals $s$.
\end{enumerate}
and the families $(C,D,z)$ and $(C',D',z')$ are said to be equivalent, if there
exists an isomorphism $C\to C'$ (over $S$) such that
the first family goes to the second under the naturally induced
morphisms.
\end{defn}

By \cite{MP} it is known that the so called ``Krichever map'' is in
fact the following morphism of functors:
$$\aligned
K:\cur &\longrightarrow \gr(k((z)),k[[z]])\\
(C,D,z)&\longmapsto \limil{n}\pi_*\o_C(n)
\endaligned$$
It is also known that $K$ is an immersion and that there exists a
locally closed subscheme of $\grv$ representing $\cur$ (which we will
denote again by $\cur$).

Let us recall another construction very similar to the above one.  Set
$m=(C,p,z)\in\cur(\spk)$, and consider the functor:
$$S\rightsquigarrow \tilde\pic(C,p)=
\left\{(L,\phi)\,\vert\,\gathered  L\in\pic(C)^{\bullet}(S)
\text{ and $\phi$ is a} \\
\text{formal trivialization of $L$ around $p$}
\endgathered\right\}$$
Define the morphism:
$$\aligned
K_m:\tilde\pic(C,p) &\longrightarrow \gr(k((z)),k[[z]])\\
(L,\phi)&\longmapsto \limil{n}\pi_* L(n)
\endaligned$$
which is also usually called ``Krichever functor''. (For a more
detailed study of $\tilde\pic(C,p)$ and $K_m$ see \cite{Al}).

Note that while $K$ is very well adapted to study of the moduli
space of curves (\cite{MP,N}), the other one, $K_m$, is good for the
study of Jacobian varieties and their subvarieties (\cite{Sh}). See
also \cite{Mul2}.

\subsection{The Formal Group $\Gamma$}\label{subsec:gamma}
Let us now recall some basic facts about the formal group $\Gamma$ (for
a complete study and definitions see \cite{AMP}). $\Gamma$ is
defined as the formal group scheme $\Gamma_-\times{\Bbb G}_m\times
\Gamma_+$ over
$\spk$, where $\Gamma_-$ is the formal scheme representing the functor
on groups:
$$S\rightsquigarrow \Gamma_-(S)=\left\{\gathered
\text{ series }\,a_n\,z^{-n}+\dots+a_1\,z^{-1}+1\\
\text{ where }a_i\in H^0(S,\o_S)\text{ are }\\
\text{ nilpotents and $n$ is arbitrary }
\endgathered\right\}$$
${\Bbb G}_m$ is the multiplicative group, and the scheme $\Gamma_+$
represents:
$$S\rightsquigarrow \Gamma_+(S)=\left\{\gathered
\text{ series }1+a_1\,z+a_2\,z^2+\dots\\
\text{ where }a_i\in H^0(S,\o_S)
\endgathered\right\}$$
The group laws of $\Gamma_-$ and $\Gamma_+$ are those induced by the
multiplication of series. Note also that there exists a natural inclusion
of $\Gamma$ in the identity connected component of:
$$S\rightsquigarrow
H^0(S,\o_S)((z))^\ast\,\overset{\text{\tiny\rm
def}}{=}\,H^0(S,\o_S)[[z]][z^{-1}]^*$$
which is an isomorphism when $char(k)=0$.

Further, $\Gamma_-$ is the inductive limit of the schemes:
$$S\rightsquigarrow \Gamma^n_-(S)=\left\{\gathered
\text{ series }\,a_n\,z^{-n}+\dots+a_1\,z^{-1}+1\\
\text{ where }a_i\in H^0(S,\o_S)\text{ and the $n^{\text{th}}$}\\
\text{ power of the ideal $(a_1,\dots,a_n)$ is zero }
\endgathered\right\}$$
in the category of formal schemes.

Observe now that there exist two actions of $g(z)\in\Gamma$ in $V$;
namely, the one given by homotheties:
$$\aligned H_g:V& \to V\\h(z)&\mapsto g(z)\cdot h(z)\endaligned$$
and the one defined by the automorphism of $k$-algebras:
$$\aligned U_g:V& \to V\\z&\mapsto z\cdot g(z)\endaligned$$

\begin{rem}\label{rem:sigma}
It is known that $U:g\mapsto U_g$ establishes a bijection
$k[[z]]^*=\Gamma_+(k)\iso\operatorname{Aut}_{\text{$k$-alg}}(k((z)))$.
Recall from \cite{Bo} (chapter III, \S4.4) that given a $k$-algebra
automorphism $\sigma$ of $k((z))$, there exists a unique
$g(z)\in\Gamma_+(k)$ such that
$U_g\circ \sigma=\sigma_0$ (where $\sigma_0$ is the $k$-algebra
automorphism  of $k((z))$ given by $z\mapsto -z$; that is, it is
possible to ``normalize'' $\sigma$ such that $\sigma^*(g(z))=g(-z)$.
\end{rem}

\begin{rem}
Now, $g\in\Gamma_+$ acts on $\tilde\pic(C,p)$ sending $(L,\phi)$ to
$(L,H_g\circ \phi)$. And hence the projection morphism:
$$\tilde\pic(C,p) @>p_1>> \pic(C)$$
may be interpreted as a principal bundle of group $\Gamma_+$. Now
comparing the zero locus of sections of $\detd_V$ and $\o(\Theta)$, one
deduces:
$$\detd_V\vert_{\tilde\pic(C,p)}\,\iso\, p_1^*\o_{\pic(C)}(\Theta)$$
which allows one to write the $\tau$-function of the point $U=K(C,p,z)$
(restricted to $\tilde\pic(C,p)$ via $K$) in terms of the theta
function of the Jacobian of $C$.

For explicit formulas relating $\tau$-functions and theta functions of
Riemann surfaces, see \cite{Kr,Sh} (see also \cite{Sh2} for the case of
Pryms).
\end{rem}

\begin{rem}
For other constructions and properties of the group $\Gamma$ see
\cite{AMP,C,KSU,PS,SW}.
\end{rem}

\section{Formal Prym variety and BKP hierarchy}

Observe now that $\Gamma$, and hence all the above-mentioned
subgroups,  acts on $\grv$ by homotheties and that on the set of
rational points ${\Bbb G}_m$ acts trivially, and
$\Gamma_+$ freely. Recall from \cite{AMP,KSU} that $\Gamma$ behaves
like the Jacobian of the formal curve. Our goal is then to achieve an
anologous result for the case of Pryms. This arises from the answer of
the following question: which is the maximal subgroup of $\Gamma$
acting on
$\gr^I_\sigma(V)$?.

\begin{thm}
The maximal subgroup of $\Gamma$ acting on $\gr^I_\sigma(V)$ is:
$$\Pi_\sigma\,=\,\{g(z)\in\Gamma\,\vert\, g(z)\cdot\sigma^*g(z)=1\,\}$$
which is a subscheme of $\Gamma$.
\end{thm}

\begin{pf}
Observe that the homothety by $g(z)\in\Gamma$ restricts to a
automorphism of $\gr^I_\sigma(V)$ if and only if:
$$g(z)\cdot U\in\gr^I_\sigma(V)\qquad\text{for all }U\in\gr^I_\sigma(V)$$
or, what amounts to the same:
$$ g(z)\cdot U \,=\, \left(g(z)\cdot U\right)^{\perp}\qquad
\text{for all } U\in\gr^I_\sigma(V)$$
 Recalling the
definition of the metric: $(f,g)\mapsto \res_{z=0}
f(z)\cdot\sigma^*g(z)dz$, one has:
$$\left(g(z)\cdot U\right)^{\perp}\,=\, (\sigma^*g(z))^{-1}\cdot
U^{\perp}$$
Note that $U=U^{\perp}$, since $U\in\gr^I_\sigma(V)$. And one concludes that
$g(z)\cdot\sigma^*g(z)\cdot U=U$ for all $U\in\gr^I_\sigma(V)$ and
hence $g(z)\cdot\sigma^*g(z)$ must be equal to 1.
\end{pf}

\begin{defn}
The formal Prym variety is the formal group scheme:
$$\Pi^\sigma_-\,\overset{\text{\tiny{def}}}{=}\,\Pi_\sigma\cap\Gamma_-$$
\end{defn}

It is therefore natural to consider $\Pi^\sigma_-$ instead of $\Gamma_-$
in the study of $\gr^I_\sigma(V)$, and hence in the study of Pryms.
Recall that the action of $\Gamma_-$ on $\grv$ is essential in the
definition of the $\tau$-function and the Baker-Akhiezer function of a
point $U\in\grv$. But for $\gr^I_\sigma(V)$ we must restrict this action
to $\Pi^\sigma_-$. Denote by $\mu^\sigma_U$ the restriction of $\mu_U$
to
$\Pi^\sigma_-$; that is:
$$\begin{array}{ccccc}
\mu_U:\Gamma_-\times\{U\} & \hookrightarrow & \Gamma_-\times\grv
& \to & \grv \\
\phantom{xxx}\cup & & \cup\; & & \cup \\
\mu^\sigma_U:\Pi^\sigma_-\times\{U\}  &  \hookrightarrow &
\Pi^\sigma_-\times\gr^I_\sigma(V) & \to & \gr^I_\sigma(V)
\end{array}$$
These actions are the cornerstone of \S4, where they will be studied at
the tangent space level.

\begin{defn}
The $\bar\tau$-function of a point $U\in\gr^I_\sigma(V)$ is the section
$({\mu^\sigma_U})^*\bar\Omega_+$ of $({\mu^\sigma_U})^*\Pf$.
\end{defn}

\begin{thm}\label{thm:taus}
$$\tau_U\vert_{\Pi^\sigma_-}\,=\,\lambda\cdot \bar\tau_U^2\qquad
\lambda\in k^*$$
\end{thm}

It is known that the KP hierarchy is a system of partial differential
equations for the $\tau$-function of a point $U\in\check{\mathbb
P}\Omega$. These are in fact equivalent to the Pl\"ucker equations for
the coordinates of $U$. It is thus quite natural to give the following:

\begin{defn}\label{defn:BKP}\hfill
\begin{itemize}
\item {\rm $char(k)$ arbitrary:} The BKP hierarchy is the set of
algebraic equations  defining $\gr^I_\sigma(V)$ inside $\check{\mathbb
P}\Omega$; in particular it gives,
\item {\rm $char(k)=0$:} The BKP hierarchy is the system of partial
differential equations that characterizes when a function is a
$\bar\tau$-function of a point of $\gr^I_\sigma(V)$.
\end{itemize}
\end{defn}

The relationship between the BKP hierarchy and Pryms will be clear in
{\ref{thm:shiota-prym}}.

\begin{rem}
Let us relate all the above claims to the classical results when
$char(k)=0$. Classically, the BKP hierarchy is introduced as the
system of equations obtained from the KP system making $t_i=0$ for all
even $i$.

Take the formal trivilization around $p$
equal to $\sigma_0$; that is: $\sigma_0(g(z))=g(-z)$ for all
$g(z)\in\Gamma$. Then, using the isomorphism of $\Gamma$ with an
additive group given by the exponential map, one has that the set of
$A$-valued points of
$\Pi^0_-$ is (we write only $0$ instead of $\sigma_0$):
$$\Pi^0_-(\sp(A))\,=\,\left\{\gathered
\text{ series }\exp\big(\sum^{n} \Sb i=1 \\ \text{$i$ odd} \endSb
a_iz^{-i}\big)  \text{ where $a_i\in A$}\\
\text{ is nilpotent and $n>0$ arbitrary}\endgathered\right\}$$

  Some well known results, such as formula 1.9.8 of \cite{DKJM}, are now
a consequence of the relationship between the
$\tau$-function of $U\in\gr^I_0(V)$ as a point of $\grv$ and the
$\bar\tau$-function as a point of $\gr^I_0(V)$ given in Theorem
{\ref{thm:taus}}.

These connections of KP and BKP, of $\grv$ and $\gr^I_0(V)$, of
$\detd_V$ and $\Pf$, and of $\tau$ and $\bar\tau$ (given above) justify
the expression given in \cite{DKJM2} of a tau-function for the BKP in
terms of theta functions of a Prym variety\footnote{It
is known that the restriction of a theta function of a Jacobian variety
is the square of a theta function of the Prym when the involution
has two fixed points (\cite{Mum}).}, as well as the methods of \cite{H}
and \cite{O} based on pfaffians of matrices in order to construct
solutions for the BKP.
\end{rem}

\section{Geometric Characterizations}

Geometric characterizations of Jacobians and Pryms offered in several
papers (\cite{Mul1,Sh,LiMu,Sh2}) are based on the study of an action of
a group on a space at the tangent space level and are therefore
suitable for being expressed in terms of dynamical systems. Roughly,
the group is a subgroup of the linear group and the space is the
Grassmannian (or a quotient of it), and the way to relate Jacobians and
Pryms with the Grassmannian is through the Krichever functor.

Since we aim to give a scheme-theoretic generalization of Theorem 6 of
\cite{Sh}, \S2.5 of \cite{Sh2} and Theorem 5.14 of \cite{LiMu}, we
should use only algebraic methods. In Shiota's paper, the action of
the group is given by analytic techniques, while the
characterization of Pryms given by Li and Mulase involves quotients  of
the Grassmannian that do not need to be algebraic.

Our approach therefore needs  to use the notion that the action of
$\Gamma$ in
$\grv$ ($\Pi_-^\sigma$ on $\gr_\sigma^I(V)$) is algebraic and that
there is no need to use quotient spaces because of Lemma
{\ref{lem:orbit}}. Essentially, this Lemma implies that the
orbit of a point of $\grv/\Gamma$ under $\Gamma$ coincides with
that of a preimage in $\grv$ under $\Gamma_-$. Although the methods are
algebraic, it seems quite natural to use the language of
dynamical systems when working at the tangent space level.

\subsection{Dynamical systems and the Grassmannian}

Observe that the action:
$$\aligned \Gamma\times\grv& @>\mu>> \grv \\
(g,U)&\mapsto g\cdot U\endaligned$$
canonically induces  a system of partial differential equations
(p.d.e.) on $\grv$. Taking $\sp(k[\epsilon]/(\epsilon^2))$-valued
points, and using the canonical identification:
$$T\grv\,\iso\,\hom(\L,\hat V/\L)$$
we obtain a morphism of functors on groups:
$$\aligned T_{\{1\}}\Gamma & @>d\mu>> \hom(\L,\hat V/\L) \\
1+\epsilon g &\longmapsto (\L\hookrightarrow\hat  V @>\cdot g>>
\hat V\to\hat  V/\L) \endaligned$$
Denote by $\mu_-$ ($d\mu_-$) the restriction
$\mu\vert_{\Gamma_-}$ ($d\mu\vert_{\Gamma_-}$ respectively).

Moreover, note that the map $g$ to $1+\epsilon g$ gives an isomorphism
of functors $\hat V\simeq T_{\{1\}}\Gamma$, where $\hat V(S)=\lim
(V/A\underset{k}\otimes\o_S)$. Also, the kernel of $d\mu_-$ at a point
$U$ is the maximal sub-$k$-algebra of $V$ acting (by homotheties) on $U$.

\begin{defn}
Given a subbundle $E\subseteq T\grv$, a finite dimensional solution of
the p.d.e. associated with $E$ at a point $U$
will be a finite dimensional subscheme $X\subseteq\grv$ containing
$U$ and such that $E_U\simeq T_{U}X$.
\end{defn}

\begin{rem}
It is convenient to consider not only finite
dimensional subschemes as solutions but also algebraizable formal
schemes.
 \end{rem}

\subsection{Characterization of Jacobian varieties}

The goal of this subsection is to prove the following generalization of
the Theorem 6 of \cite{Sh}:

\begin{thm}\label{thm:shiota}
A necessary and sufficient condition for a  rational
point $U\in\grv$ to lie in the image of the Krichever map $K_m$ (for a
point $m\in\cur$) is that there exists a finite dimensional solution of
the p.d.e. $\im{d\mu_-}$ at the point $U$.
\end{thm}

Proof of the theorem is a direct consequence of the following two
lemmas, that are quite akin to Mulase's and Shiota's ideas
(\cite{Mul1,Sh}).

\begin{lem}\label{lem:orbit}
Let $U$ be rational point of $\grv$, and let $G(U)$ denote the
orbit of $U$ under the action of a group $G$. Then:
$$\Gamma(U)\simeq \Gamma_-(U)\times\Gamma_+$$
\end{lem}

\begin{pf}
Since $\Gamma=\limi(\Gamma^n_-\times{\Bbb G}_m\times\Gamma_+)$ (as
formal schemes), it is then enough to show that the natural map:
$$\aligned
\Gamma^n_-(U)\times\Gamma_+&\to
(\Gamma^n_-\times{\Bbb G}_m\times\Gamma_+)(U)\\
\left(f_-U,f_+\right)&\mapsto (f_-\cdot f_+)U
\endaligned$$
is an isomorphism; or, what amounts to the same: if $f_-U=f_+U$,
then $f_-=f_+=1$. Since $U$ is a rational point, there
exists an element $g=\sum_{i\ge m}c_i\,z^i\in U$ ($c_m\ne 0$) with
$m$ maximal (this $m$ will be called order of $g$); note also that we
can assume $c_m=1$. Note that homotheties map elements of maximal order
into elements of maximal order, and therefore  $f_-\cdot g=f_+\cdot g$,
from which one deduces $f_-=f_+$, and hence
$f_-,f_+\in\Gamma_-\cap\Gamma_+=\{1\}$.

Finally, the natural structure of the formal scheme of
$\Gamma_-(U)$ is equal to $\cup_{n>0}\Gamma^n_-(U)$
(where $\Gamma^n_-(U)$ denotes the schematic image of
$\Gamma^n_-\times\{U\}\to\grv$).

Let $\Gamma^n_-$ be
$\sp\left(k[x_1,\dots,x_n]/(x_1,\dots,x_n)^n\right)$, $\o$ the structural
sheaf $\o_{\grv}$, $I_U^n$ the kernel of  $\o\to
k[x_1,\dots,x_n]/(x_1,\dots,x_n)^n$, and $u^{n'}_n$ the morphism
$\o/I_U^{n'}\to\o/I_U^n$ ($n'>n$).

In order to show that $\Gamma_-(U)=\sf(A)$ ($A$ being $\limp\o/I_U^n$),
it remains only to check that (\cite{EGA} I.10.6.3) $u^{n'}_n$ is
surjective and $\ker(u^{n'}_n)$ is nilpotent. However, both are obvious.

Note that as a by-product we have that the topology of $A$
is given by the ideals $J_n=\ker(A\to \o/{I_U^n})$ and
that the definition ideal is $J=\limpl{n}(x_1,\dots,x_n)$.
\end{pf}

\begin{lem}\label{lem:five}
Let $U$ be a rational point of $\grv$. Then the following conditions
are equivalent:
\begin{enumerate}
\item $\Gamma_-(U)$ is algebraizable,
\item $\dim_k(J/J^2)<\infty$,
\item $\dim_k(T_U\Gamma_-(U))<\infty$,
\item $\dim_k\im(d\mu_-)<\infty$,
\item there exists a rational point $m=(C,p,z)$ of $\cur$ and a
pair $(L,\phi)$ of $\tilde\pic(C,p)(\spk)$ such that $K_m(L,\phi)=U$.
\end{enumerate}
\end{lem}

\begin{pf}
$1\implies 2$: Recall that
algebraizable (\cite{H} II.9.3.2) means that the formal scheme is
isomorphic to the completion of a noetherian scheme along a closed
subscheme, but the completion of a noetherian ring with respect to an
ideal is noetherian (\cite{At} 10.26); hence $A$ is noetherian.

Recall now from \cite{EGA} 0.7.2.6 that if $A$ is an
admissible linearly topologized ring, and $J$ a definition ideal,  $A$ is
noetherian if and only if $A/J$ is so and
$J/J^2$ is a finite type $A/J$-module.

Since $A/J\simeq k$, we conclude.

$2\implies 3$: Note that
$$\aligned T_U\Gamma_-(U)&=
\hom_{\text{for-esq}}(\sp(k[\epsilon]/(\epsilon^2),\sf(A))=\\
&=\hom \Sb \text{topological} \\ \text{$k$-algebras} \endSb
(A,k[\epsilon]/(\epsilon^2)) \subset
\hom_{\text{$k$-algebras}}(A,k[\epsilon]/(\epsilon^2))\endaligned$$
which is isomorphic to $\left( J/J^2\right)^*$.

$3 \iff 4$ By the very definition, the morphism $\Gamma_-\to\grv$
factorizes through $\Gamma_-(U)$, and therefore:
$$\dim_k\im(d\mu_-)\,\leq\,\dim_k T_U\Gamma_-(U)<\infty$$

$4\implies 5$: From 4 we have $\dim_k\ker(d\mu_-)=\infty$. Note,
moreover, that $B=\ker(d\mu_-)\subseteq V$ is a integral $k$-algebra of
transcendence degree 1, since $B_{(0)}= V$ such that $U$ is a free
$B$-module of rank 1. Standard results (see \cite{SW}) show
how from the pair $(B,U)$ one can construct the data
$(C,p,z)\in\P(\spk)$ and
$(L,\phi)\in\tilde\pic(C,p)(\spk)$ such that $K_m(L,\phi)=U$.

$5\implies 4$: Easy.

$3\implies 1$: Let $\sp(A_n)$ be the schematic image of
$$\Gamma_-^n=\sp(k[x_1,\dots,x_n]/(x_1,\dots,x_n)^n)\to\grv$$
 By the
properties of formal schemes, one has:
$$ T_U\Gamma_-(U)\quad=\quad \bigcup_{n>0} T_U \sp(A_n)$$
And therefore, if $\dim_k (T_U\Gamma_-(U))=d$ then $\dim_k (T_U
\sp(A_n))=d$ for all $n>>0$.

Denote with $J_{(n)}$ the maximal ideal of
$A_n$, since $T_U \sp(A_n)\simeq (J_{(n)}/J_{(n)}^2)^*$ and $J=\limp
J_{(n)}@>\pi_n>> J_{(n)}$ is surjective, there exist elements $\bar
y_1,\dots,\bar y_d\in J$ such that $<\{\pi_n(\bar y_1),\dots,
\pi_n(\bar y_d)\}>=J_{(n)}/J_{(n)}^2$. By Nakayama's lemma one has an
epimorphism:
$$\aligned p_n:k[y_1,\dots,y_d]&\to A_n\\
y_i\quad&\mapsto\pi_n(\bar y_i)\endaligned$$

It is now straightforward to see that $\{p_n\}$ is compatible with the
natural  morphism $A_m\to A_n$ for $m>n$. One concludes that $A_n\simeq
k[y_1,\dots,y_d]/I_n$ ($I_n$ being $\ker p_n$) and that
$I_m\subset I_n$ for $m> n$.

Recall now that $A=\limp A_n$ and thus $A\simeq k[[y_1,\dots,y_d]]$;
the topology on $A$ induced by $\ker(A\to A_n)$ coincides with the
$(y_1,\dots,y_d)$-adic; that is, the formal scheme $\sf(A)$ is
algebraizable.
\end{pf}

\begin{rem}
Assume that the conditions of Lemma {\ref{lem:five}} are satisfied for a
point $U$. Then the orbit $\Gamma_-(U)$, which is a formal scheme,
 is the desired finite dimensional solution. Let us give another
interpretation. Let $m=(C,p,z)\in\cur$ and $(L,\phi)\in\tilde\pic(C,p)$
such that $K_m(L,\phi)=U$. It is straightforward to check that the
p.d.e. defined by $\im(d\mu_-)$ is in fact the KP hierarchy and hence
$\tilde\pic(C,p)$ is a finite dimensional solution of the KP hierarchy,
modulo the action of $\Gamma_+$.
\end{rem}

\begin{rem}
Let $char(k)=0$, $U\in\grv$ satisfying one condition of the
previous lemma and $(C,p,z,L,\phi)$ that given by the fifth condition.

Then, the morphism
$T_{\{1\}}\Gamma^n_-\to T_U\grv$ is essentially:
$$H^0(C,\o_C(np)/\o_C)\to H^1(C,\o_C)\simeq T_L \pic(C)$$
which is deduced from the cohomology sequence of:
$$0\to\o_C\to\o_C(np)\to \o_C(np)/\o_C\to 0$$
and from the proof of the Lemma it is easy to conclude that
$\Gamma_-(U)$ is canonically isomorphic to the formal completion of the
Jacobian of $C$ along $L$.

When $k=\Bbb C$ and $C$ is a smooth curve,  we can interpret the
algebraic variety $\pic(C)$ as a compact Lie group. Now using the
exponential map, it is not hard to see that our condition of ``finite
dimensional solution of the p.d.e.'' is
equivalent to Shiota's ``compact solution for the KP hierarchy'' (see
Theorem 6 of \cite{Sh}), and that the  morphism
$T_{\{1\}}\Gamma^n_-\to T_U\grv$ is the one studied in depth in its
\S2, especially in Lemmas 2 and 4 of \cite{Sh}.
\end{rem}

\begin{rem}
Note further that the point $m=(C,p,z)\in\cur$ obtained by Lemma
{\ref{lem:five}} is not unique. However, it has a characterizing
property. From the construction, it is easily seen that the ring $B$ is
the maximal subring of $\hat V$ such that $B\cdot U= U$. This implies
that there is no morphism  $f:C'\to C$ and line bundle $L'$ such that
$f_*L'=L$ (unless $f$ is an isomorphism). See \cite{Mul2} for more
details.
\end{rem}

\subsection{Characterization of Prym varieties}
We shall say that a point $(C,p,z)\in\cur$ admits the automorphism
$\sigma$ if $\sigma:k((z))\iso k((z))$ (as a $k$-algebra) restricts to
$K(C,p,z)\iso K(C,p,z)$. Note that in this case $\sigma$ induces an
automophism of $C$ with $p$ as a fixed point (it will also be denoted by
$\sigma$).

\begin{defn}\label{defn:p0}
 Let $\P_\sigma$ denote the subfunctor of $\cur$ consisting of the data
that admit the automorphism $\sigma$. (In the next section we shall see
that it is in fact a subscheme).
\end{defn}

In this setting, the Prym variety associated with the data
$m=(C,p,z)\in\cur$ that admits an involution $\sigma$ is a subscheme of
$\tilde\pic(C,p)$ whose rational points are:
$$\tilde{\operatorname{Prym}}(C,p,\sigma)=\left\{
(L,\phi)\in\tilde\pic(C,p)\,\vert\,
\sigma^*(L)\iso
\omega_C\otimes L^{-1}\right\}$$

In this subsection we shall restrict ourselves to the situation
addressed in Example {\ref{exam:metric}}; that is:
$\sigma=\sigma_0$. Note, however, that this can always be achieved (see
Remark {\ref{rem:sigma}}). Thus, we shall assume here that
$\sigma_0^*(g(z))=g(-z))$ for all $g(z)\in\Gamma$, and shall remove the
super/sub-script $\sigma$ in the notations. From the above discussion
and recalling {\ref{defn:gri}}, one has the following cartesian diagram:
\begin{eqnarray*}
\tilde\pic(C,p)\phantom{xx} &
\overset{\text{\scriptsize $K_m$}}\hookrightarrow & \,\grv
\\
\cup\phantom{xxxx} & & \phantom{xx} \cup \\
\tilde{\operatorname{Prym}}(C,p,\sigma_0)
& \overset{\text{\scriptsize $K_m$}}\hookrightarrow & \gr^I_0(V)
\end{eqnarray*}

Let $\mu^0_-$ be the restriction of $\mu_-$ to $\Pi^0_-$, and
let $d\mu^0_-$ be that induced in the tangent spaces. Our version of
the Theorem 5.14 of \cite{LiMu} is the following:

\begin{thm}\label{thm:shiota-prym}
A necessary and sufficient condition for a  rational
point $U\in\gr^I_0(V)$ to lie in the image of the Krichever map
$K_m$ (for a point of $m\in\P_0$) is that there exists a finite
dimensional solution of the p.d.e. $\im{d\mu^0_-}$ at the point $U$.
\end{thm}

\begin{pf}
Observe that $\dim_k\im{d\mu^0_-}<\infty$ if and only if
$\dim_k\im{d\mu_-}<\infty$, and therefore that it is only necessary to
show that in the last condition of Lemma {\ref{lem:five}}, the
constructed data $(C,p,z)\in\cur$ admits the involution given by
$z\mapsto -z$.

First, note that since $U$ is m.t.i. we have: $<f,u>=0\quad \forall
\,u\in U\implies f\in U$. Now for an element $f\in\ker(d\mu^0_-)$ we have
$f\cdot U\subseteq U$ and therefore $<f\cdot u,v>=0$ for all $u,v\in U$.

We want to see that $f(-z)\in\ker(d\mu^0_-)$; that is,
$<f(-z)u(z),v(z)>=0$ for all $u,v\in U$. Note, however, that:
$$<f(-z)u(z),v(z)>=<u(z),f(z)v(z)>=-<f(z)v(z),u(z)>$$
and hence $f(-z)\in\ker(d\mu^0_-)$, as desired.
\end{pf}

\begin{rem}
Analogously to the case of Jacobian varieties, one has that
$\im(d\mu^0_-)$ is equivalent to the BKP hierarchy and therefore
$\Pi^0_-(U)$ is a finite dimensional solution for the BKP hierarchy. As
before, one can say that $\tilde{\operatorname{Prym}}(C,p,z)$ (modulo
the action of $\Gamma_+$) is a finite dimensional solution too.
\end{rem}

\begin{rem}
Note that in the proof of Lemma {\ref{lem:five}} a $k$-algebra $B$ is
constructed that turns out to be the ring $H^0(C-p,\o_C)$. If we now
assume that $U\in\gr^I_0(V)$, $B$ has an involution. It is easy to
check that in this situation $B\in\gr_0(k((z)),k[[z]])$. Since we know
the structure of this Grassmannian (see Example {\ref{exam:metric}}),
the projection on the first factor is precisely:
$$B'\,=\,B\cap k((z^2))\,\in\,\gr_0(k((z^2)),k[[z^2]])$$

Now, the curve $C'$ constructed from $B'$ is the quotient of $C$ with
respect to the involution $\sigma$. Note that $U$ is a rank 2
free $B'$-module, and that the sheaf induced by $U$ over $C'$ is the
direct image of $L$ by $C\to C'$.
\end{rem}

\section{Equations for the Moduli Space of Prym varieties}

Although the notion of Prym variety is more general, here we shall
restrict our study to those coming from a curve and an involution with
at least one fixed point.

\begin{defn}
The functor of Prym varieties, $\P$, is the sheafication of the
following functor on the category of $k$-schemes:
$$S\rightsquigarrow\left\{(C,D,z,\sigma)\,\vert\,
\gathered (C,D,z)\in\cur(S),\text{ $\sigma$ is an involution} \\
\text{that induces an automorphism of} \\
\text{the formal completion of $C$ along $D$}
\endgathered
\right\}$$
(up to isomorphisms).
\end{defn}

Given $(C,D,z,\sigma)\in\P(\spk)$, note that $\sigma$ induces an
automorphism of $k((z))$, and that the action of $\Gamma_+$ on the group
$\operatorname{Aut}_{\text{$k$-alg}}(k((z)))$ (via $U$) is transitive
and free. One  therefore has a bijection (set-theoretic)
$\P\simeq\P_0\times\Gamma_+$.

From Example {\ref{exam:metric}} and Definition {\ref{defn:p0}} one
now has the following easy but fundamental result:

\begin{thm}\label{thm:char-p0}
Via the Krichever map, one has:
$$\P_0\simeq \cur \underset{\grv}\times\gr_0(V)$$
\end{thm}

\begin{cor}
The functor $\P_0$ is representable by a locally closed subscheme of the
infinite Grassmannian $\grv$; namely, that whose ($S$-valued) points $U$
 satisfy:
\begin{enumerate}
\item $\o_S\subset U$,
\item $U\cdot U\subseteq U$ ($\cdot$ being the product of $\o_S((z))$),
\item the map $\o_S((z))\to \o_S((z))$ defined by $z\mapsto -z$
restricts to an isomorphism $U\iso U$.
\end{enumerate}
\end{cor}

\begin{pf}
Recall from \cite{Al,MP} that the first two conditions are locally
closed. The third condition is closed since it is where the identity and
the involution of $\grv$ given by $z\mapsto -z$ coincide; recall also
that
$\grv$ is separated.
\end{pf}

We can now state a theorem characterising the points of $\gr_0(V)$ in
terms of bilinear identities; this is an analogous result to the
characterization of $\grv$ of \cite{DKJM,F,MP}.

\begin{thm}[Bilinear Identities]\label{thm:bil-ident}
$$U\in\gr_0(V)\quad\iff\quad\left\{
\gathered \res_{z=0}\psi_U(z,t)\psi^*_U(z,t')\frac{dz}{z^2}=0 \\
\res_{z=0}\psi_U(-z,t)\psi^*_U(z,t')\frac{dz}{z^2}=0 \endgathered\right.
\,\text{ for all $t,t'$}$$
\end{thm}

\begin{pf}
Note that the third condition in the corollary is equivalent (for
the rational points) to saying that: $\psi_U(z,t)=-\psi_U(-z,t)$ since
for a point $U\in\grv$ one has that the Baker-Akhiezer function,
$\psi_U(z,t)$, may be written as a series of the form
$z\cdot\sum_{i>0}\psi_U^{(i)}(z)p_i(t)$ where $p_i(t)$ are universal
polynomials (they do not depend on $U$) and $\{\psi_U^{(i)}(z)\}_{i>0}$
is a basis of $U$ as a $k$-vector space (see \cite{MP}).

Recalling the property
$\res_{z=0}\psi_U(z,t)\psi^*_{U'}(z,t')\frac{dz}{z^2}=0$ if and only if
$U=U'$, one concludes.
\end{pf}

This result holds when $char(k)\neq 2$; however, it can be translated
in the language of differential equations when $char(k)=0$. We arrive
at two sets of differential equations in terms of $\tau$-functions
defining the set of rational points of $\gr_0(V)$ and
$\P_0$.

For a Young diagram $\lambda$, denote by $\chi_\lambda$ its associated
Schur polynomial. If the diagram has only one row of length $\beta$, then
denote the corresponding Schur polynomial by $p_\beta$. Let $t$ be
$(t_1,t_2,\dots)$ and $\tilde\partial_t$ be
$(\partial_{t_1},\frac12\partial_{t_2},\frac13\partial_{t_3},\dots)$.
Let $D_{\lambda,\alpha}$ be $\sum_\mu\chi(\tilde\partial_t)$ where
the sum is taken over the set of Young diagrams $\mu$ such that
$\lambda-\mu$ is a horizontal $\alpha$-strip.

\begin{thm}[P.D.E. for $\gr_0(V)$]\label{thm:pde-gr0}
A function $\tau(t)$ is the $\tau$-function of a rational point
$U\in\gr_0(V)$ if and only if it satisfies the following infinite set
of differential equations (indexed by a pair of Young diagrams
$\lambda_1,\lambda_2$):
$$\left(\sum
p_{\beta_1}(-\tilde\partial_t)D_{\lambda_1,\alpha_1}(\tilde\partial_t)\cdot
p_{\beta_2}(-\tilde\partial_{t'})D_{\lambda_2,\alpha_2}(\tilde\partial_{t'})
\right)\vert_{t=t'=0}\tau_U(t)\cdot \tau_U(t')\,=\,0$$
where the sum is taken over the 4-tuples
$\{\alpha_1,\beta_1,\alpha_2,\beta_2\}$ of integers such that
$-\alpha_1+\beta_1-\alpha_2+\beta_2=1$, $-\alpha_1+\beta_1$ is even.
\end{thm}

\begin{pf}
Recall from \cite{MP} how the KP hierarchy is deduced from the Residue
Bilinear Identity. Apply the same procedure to the identities in the
preceding Theorem, and add and subtract both identities.
\end{pf}

We finish with the partial differential equations for
$\tau$-functions that characterize $\P_0$ as a subscheme of
$\check{\mathbb P}\Omega$.

\begin{thm}[P.D.E. for $\P_0$]\label{thm:pde-p0}
A function $\tau(t)$ is the $\tau$-function of a point $U\in\P_0$ if
and only if it satisfies the following infinite set of differential
equations:
\begin{enumerate}
\item the p.d.e. of Theorem {\ref{thm:pde-gr0}},
\item $$P(\lambda_1,\lambda_2,\lambda_3)\vert \Sb t=0 \\ t'=0 \\ t''=0 \endSb
\left(\tau_U(t)\cdot \tau_U(t')\cdot\tau_U(t'')\right)\,=\,0$$
for every three Young diagrams $\lambda_1,\lambda_2,\lambda_3$, where
$P(\lambda_1,\lambda_2,\lambda_3)$ is the differential operator defined
by:
$$\sum
p_{\beta_1}(\tilde\partial_t)D_{\lambda_1,\alpha_1}(-\tilde\partial_t)\cdot
p_{\beta_2}(\tilde\partial_{t'})D_{\lambda_2,\alpha_2}(-\tilde\partial_{t'})
\cdot
p_{\beta_3}(-\tilde\partial_{t''})D_{\lambda_3,\alpha_3}(\tilde\partial_{t''})
$$
where the sum is taken over the 6-tuples
$\{\alpha_1,\beta_1,\alpha_2,\beta_2,\alpha_3,\beta_3\}$ of integers such that
$-\alpha_1+\beta_1-\alpha_2+\beta_2-\alpha_3+\beta_3=2$,
$-\alpha_1+\beta_1$ is even,
\item the p.d.e.'s:
$$\left(\sum_{\alpha\geq 0}
p_{\alpha}(-\tilde\partial_t)D_{\lambda,\alpha}(\tilde\partial_t)
\right)\vert_{t=0}\tau_U(t)=0\quad
\text{for all Young diagram }\lambda$$
\end{enumerate}
\end{thm}

\begin{pf}
By Theorems {\ref{thm:char-p0}} and {\ref{thm:pde-gr0}}, it is enough to
recall from \cite{MP} the p.d.e. defining $\cur$ in the infinite
Grassmannian.
\end{pf}

Note that a theta function of the Prym variety associated to a 4-uple
$(C,p,z,\sigma)$ (where $\sigma$ is the given by $z\mapsto -z$) satisfy
these differential equations which are not a consequence of the BKP
hierarchy.



\begin{thebibliography}{XXX}

\bibitem[Al]{Al} A \'Alvarez-V\'azquez, {\it Estructuras
aritm\'eticas de las curvas algebraicas}, Ph.D. Thesis, Salamanca, June
1996

\bibitem[AM]{At} M Atiyah, I G MacDonald {\it Introduction to
commutative algebra}, Addison Wesley Publ. Company, London (1969)

\bibitem[AMP]{AMP} A \'Alvarez V\'azquez; J M Mu\~noz Porras and F J
Plaza Mart\'{\i}n, {\it The algebraic formalism of soliton equations
over arbitrary base fields}, alg-geom/9606009 (To appear in
``Proceedings of Workshop on Abelian Varieties and Theta Functions'',
edited by Aportaciones Matem\'aticas de la Sociedad Matem\'atica
Mexicana)

\bibitem[B]{B} D Borthwick, {\it The Pfaffian line bundle},  Comm.
Math. Physics  {\bf 149} (1992), pp. 463--493

\bibitem[Bo]{Bo} N Bourbaki,{\it Commutative Algebra},
Hermann (1972)

\bibitem[BS]{BS} A A Beilinson and V V Schechtman, {\it
Determinant Bundles and Virasoro Algebras} Commun. Math. Phys. {\bf
118} (1988), pp. 651-701

\bibitem[C]{C} C Contou-Carr\'ere, {\it Jacobienne locale, groupe de
bivecteurs de Witt universel, et symbole mod\'er\'e},  C. R. Acad. Sci.
Paris, s\'erie I  {\bf 318} (1994), pp. 743-746

\bibitem[DKJM]{DKJM} E Date; M Kashiwara; M Jimbo and T Miwa, {\it
Transformation groups for soliton equations}, Proceedings of RIMS
Symposium, Kyoto 1981

\bibitem[DKJM2]{DKJM2} E Date; M Kashiwara; M Jimbo and T Miwa, {\it
Quasi-periodic solutions of the orthogonal KP equations ---
Transformation groups for soliton equations V}, Publ. RIMS, Kyoto Univ.
{\bf 18} (1982), pp.1111-1119

\bibitem[EGA]{EGA} A Grothendieck  and J A Dieudonn\'e,
{\it El\'ements de g\'eom\'etrie alg\'ebrique {\bf I}},
Springer--Verlag (1971)

\bibitem[F]{F} J D Fay, {\it Bilinear Identities for Theta
Functions},  Preprint

\bibitem[H]{H} R Hirota, {\it Solitons solutions to the BKP
equations. I. The pfaffian technique}, Journal of the Physical Society of
Japan  {\bf 58}, 7, (1989), pp. 2285-2296

\bibitem[K]{Kr} I M Krichever, {\it Methods of algebraic geometry in
the theory of non-linear equations}, Russian Math. Surveys  {\bf 32}:6,
(1977), pp. 185-213

\bibitem[KSU]{KSU} T Katsura, Y Shimizu and K Ueno, {\it Formal
Groups and Conformal Field Theory over $\Z$}, Advanced Studies in
Pure Mathematics {\bf 19} (1989), Integrable systems in quantum field
theory and statistical mechanics, pp. 347-366

\bibitem[KM]{KM} F Knudsen and D Mumford, {\it The projectivity of
the moduli space of stable curves I: preliminaries on {\it det} and
{\it div}}, Math. Scand.{\bf 39} (1976), pp. 19--55

\bibitem[LM]{LiMu} Y Li and M Mulase, {\it Category of morphisms of
algebraic curves and a characterization of Prym varieties},
Max-Planck Preprint 1992 and alg-geom/9203002

\bibitem[M1]{Mul1} M Mulase, {\it Cohomological structure in soliton
equations and jacobian varieties}, J. Differetnial Geom. {\bf 19}
(1984), pp. 403-430

\bibitem[M2]{Mul2} M Mulase, {\it A correspondence between an Infinite
Grassmannian and arbitrary vector bundles on algebraic curves}, Proc.
Symp. Puer Math. {\bf 49} Series A (1989), pp. 39-50

\bibitem[Mu]{Mum} D Mumford, {\it Prym varieties}, Contributions to
Analyis {\bf } (1974), pp. 325-350

\bibitem[MP]{MP} J M Mu\~noz Porras and F J Plaza Mart\'{\i}n, {\it
Equations of the moduli space of pointed curves in the infinite
Grassmannian} (to appear)

\bibitem[N]{N} Y Namikawa, {\it A Conformal Field Theory on Riemann
Surfaces Realized as Quantized Moduli Theory of Riemann Surfaces},
Proceedings of Symposia in Pure Mathematics {\bf 49} (1989), Part I,
pp. 413-443

\bibitem[O]{O} Y Ohta, {\it Paffian solutions for the Veselov-Novikov
equations}, Journal of the Physical Society of Japan  {\bf 61}, 11,
(1992), pp. 3928-3933

\bibitem[P]{P} F J Plaza Mart\'{\i}n, {\it Picard group and
automorphism group of $\gr(k((z)))$} (to appear)

\bibitem[PS]{PS} A Pressley and G Segal, {\it Loop Groups}, Oxford
University Press

\bibitem[S]{Sh} T Shiota, {\it Characterization of Jacobian
varieties in terms of soliton equations},  Inventiones Mathematicae
{\bf 83} (1986), pp. 333--382

\bibitem[S2]{Sh2} T Shiota, {\it Prym varieties and soliton equations},
Infinite-dimensional Lie algebras and groups,  Adv. Ser. Math. Phys.,
7, (Luminy-Marseille, 1988), pp. 407--448, World Sci. Publishing

\bibitem[SW]{SW} G Segal and G Wilson, {\it Loop groups and equations
of KdV type},  Publ. Math. I.H.E.S. {\bf 61} (1985), pp. 3--64

\bibitem[SS]{SS} M Sato and Y Sato, {\it Soliton equations as
dynamical systems on infinite Grassmann manifold}, Lecture Notes
in Num. Appl. Anal. {\bf 5} (1982), pp. 259--271

\end{thebibliography}
\end{document}